\documentclass{gtart}

\def\ifplaintex{\expandafter\ifx\csname documentclass\endcsname\relax}

\def\gtp{{\mathsurround=0pt\it $\cal G\mskip-2mu$eometry \&\ 
$\cal T\!\!$opology $\cal P\!$ublications}}  

\def\recd{{\small Received:\qua\receiveddate\ifx\reviseddate\relax
\else\qquad Revised:\qua\reviseddate\fi\par}} 

\def\lognumber#1{\def\thelognumber{#1}}
\def\volumenumber#1{\def\thevolumenumber{#1}}
\def\volumeyear#1{\def\thevolumeyear{#1}}
\def\papernumber#1{\def\thepapernumber{#1}}
\def\pagenumbers#1#2{\def\startpage{#1}\def\finishpage{#2}}
\def\published#1{\def\publishdate{#1}}

\def\received#1{\def\receiveddate{#1}}
\def\revised#1{\def\reviseddate{#1}}
\def\accepted#1{\def\accepteddate{#1}}

\long\def\asciiabstract#1{\long\def\theasciiabstract{#1}}

\let\\\par\let\thelognumber\relax\let\thevolumenumber\relax
\let\thepapernumber\relax\let\thevolumeyear\relax\let\startpage\relax
\let\finishpage\relax\let\publishdate\relax\let\receiveddate\relax
\let\reviseddate\relax\let\accepteddate\relax\let\theasciititle\relax
\let\theasciiauthors\relax
\let\theasciiabstract\relax

\let\theasciiemail\relax

\ifplaintex
\font\logobig=cmssbx10 scaled 3836
\font\logomed=cmssbx10 scaled 2557
\else
\font\logobig=cmssbx10 scaled 4200
\font\logomed=cmssbx10 scaled 2800
\fi

\long\def\makeagttitle{   
\count0=\startpage
\agt\hfill      
\hbox to 45truept{\vbox to 0pt{\vglue -13truept{\logomed A\kern -.37em{\logobig 
T}\kern -.38em G}\vss}\hss}
\break
{\small Volume \thevolumenumber\ (\thevolumeyear)
\startpage--\finishpage\nl
Published: \publishdate}

\vglue .25truein

{\parskip=0pt\leftskip 0pt plus
1fil\def\\{\par\smallskip}{\Large\bf\thetitle}\par\medskip} \vglue
0.05truein

{\parskip=0pt\leftskip 0pt plus 1fil\def\\{\par}{\sc\theauthors}
\par\medskip}%
 
\vglue 0.03truein 

{\small\leftskip 25truept\rightskip 25truept{\bf Abstract}\stdspace\theabstract

{\bf AMS Classification}\stdspace\theprimaryclass
\ifx\thesecondaryclass\relax\else; \thesecondaryclass\fi\par
{\bf Keywords}\stdspace \thekeywords\par}\vglue 7truept

}   

\ifplaintex
\font\phead=cmsl9 scaled 950
\font\pnum=cmbx10 scaled 913
\font\pfoot=cmsl9 scaled 950
\headline{\vbox to 0pt{\vskip -4.5mm\line{\small\phead\ifnum
\count0=\startpage ISSN 1472-2739 (on-line) 1472-2747 (printed)
\hfill {\pnum\folio}\else\ifodd\count0\def\\{ }%
\ifx\theshorttitle\relax\thetitle\else\theshorttitle\fi\hfill{\pnum\folio}
\else\def\\{ and }{\pnum\folio}\hfill\ifx\theshortauthors\relax\theauthors
\else\theshortauthors\fi\fi\fi}\vss}}
\footline{\vbox to 0pt{\vglue 0mm\line{\small\pfoot\ifnum\count0=\startpage
\copyright\ \gtp\hfill\else
\agt, Volume \thevolumenumber\ (\thevolumeyear)\hfill\fi}\vss}}
\else
\font=cmsl9 scaled 1050
\font\lnum=cmbx10 
\font=cmsl9 scaled 1050
\makeatletter
\def\@oddhead{{\small\lhead\ifnum\count0=\startpage ISSN 1472-2739 
(on-line) 1472-2747 (printed)\hfill {\lnum\number\count0}\else\ifodd\count0
\def\\{ }\ifx\theshorttitle\relax \thetitle \else\theshorttitle\fi\hfill
{\lnum\number\count0}\else\def\\{ and }{\lnum\number\count0}
\hfill\ifx\theshortauthors\relax 
\theauthors\else\theshortauthors\fi\fi\fi}}\def\@evenhead{\@oddhead}
\def\@oddfoot{\small\lfoot\ifnum\count0=\startpage\copyright\ \gtp\hfill\else
\agt, Volume \thevolumenumber\ (\thevolumeyear)\hfill\fi}
\def\@evenfoot{\@oddfoot}
\makeatother
\fi
\let\maketitlepage\makeagttitle
\let\makeshorttitle\maketitlepage
\let\maketitle\maketitlepage

\newwrite\gtoutfile
\long\gdef\makeheadfile{  
{\def\\{, }\def\s{ }
\immediate\openout\gtoutfile head.xxx
\immediate\write\gtoutfile{To: math@arxiv.org}
\immediate\write\gtoutfile{Subject: put OR rep NNNNN:ppppp}
\immediate\write\gtoutfile{--text follows this line--}
\immediate\write\gtoutfile{Proxy-for: \ifx\theasciiauthors\relax
\theauthors\else\theasciiauthors\fi\s<\ifx\theasciiemail\relax\theemail\else\theasciiemail\fi>}
\immediate\write\gtoutfile{\noexpand\\}
\immediate\write\gtoutfile{Authors: \ifx\theasciiauthors\relax
\theauthors\else\theasciiauthors\fi}
{\def\\{ }\immediate\write\gtoutfile{Title: \ifx\theasciititle\relax
\thetitle\else\theasciititle\fi}}
\immediate\write\gtoutfile{Subj-class: GT or SG, GR etc}
\immediate\write\gtoutfile{MSC-class: \theprimaryclass\ifx\thesecondaryclass\relax\else, \thesecondaryclass\fi}
\immediate\write\gtoutfile{Journal-ref: Algebr. Geom. Topol. \thevolumenumber\s
(\thevolumeyear) \startpage-\finishpage}
\immediate\write\gtoutfile{Comments: Published by Algebraic and
Geometric Topology at}
\immediate\write\gtoutfile{\s\s\s  http://www.maths.warwick.ac.uk/agt/AGTVol\thevolumenumber/agt-\thevolumenumber-\thepapernumber.abs.html}
\immediate\write\gtoutfile{\noexpand\\}
\immediate\write\gtoutfile{}
\ifx\theasciiabstract\relax
\immediate\write\gtoutfile{\theabstract}\else
\immediate\write\gtoutfile{\theasciiabstract}\fi
\immediate\write\gtoutfile{}
\immediate\write\gtoutfile{\noexpand\\}
\immediate\write\gtoutfile{}
\immediate\closeout\gtoutfile}}  

\def\maketitlepage{\makeagttitle\makeheadfile}
\let\makeshorttitle\maketitlepage
\let\maketitle\maketitlepage

\lognumber{26}
\volumenumber{3}
\volumeyear{2003}
\papernumber{26}
\published{24 August 2003}
\pagenumbers{777}{789}
\received{18 February 2003}
\revised{23 July 2003}
\accepted{19 August 2003}

\usepackage{amsmath, amscd, amssymb}
\usepackage{latexsym}
\def \R{\mathbb{R}}
\def \Z{\mathbb{Z}}
\def \C{\mathbb{C}}

\newtheorem{theorem}{Theorem}[section]
\newtheorem{lemma}[theorem]{Lemma}

\newtheorem{corollary}[theorem]{Corollary}

\begin{document}

\title {The Chess conjecture}

\author{Rustam Sadykov}

\address{University of Florida, Department of Mathematics,\\ 358
Little Hall, 118105, Gainesville, Fl, 32611-8105, USA}

\email{sadykov@math.ufl.edu}

\begin{abstract} We prove that the homotopy class of a Morin mapping
$f\co P^p\to Q^q$ with $p-q$ odd contains a cusp mapping. This
affirmatively solves a strengthened version of the Chess
conjecture \cite{Ch},\cite{It}. Also, in view of the Saeki-Sakuma
theorem \cite{SS} on the Hopf invariant one problem and Morin
mappings, this implies that a manifold $P^p$ with odd Euler
characteristic does not admit Morin mappings into $\R^{2k+1}$ for
$p\ge 2k+1\ne 1,3,7$.
\end{abstract}

\asciiabstract{We prove that the homotopy class of a Morin mapping
f: P^p --> Q^q with p-q odd contains a cusp mapping. This
affirmatively solves a strengthened version of the Chess conjecture
[DS Chess, A note on the classes [S_1^k(f)], Proc. Symp. Pure Math.,
40 (1983) 221-224] and [VI Arnol'd, VA Vasil'ev, VV Goryunov, OV
Lyashenko, Dynamical systems VI. Singularities, local and global
theory, Encyclopedia of Mathematical Sciences - Vol. 6 (Springer,
Berlin, 1993)].  Also, in view of the Saeki-Sakuma theorem [O Saeki, K
Sakuma, Maps with only Morin singularities and the Hopf invariant one
problem, Math. Proc. Camb. Phil. Soc. 124 (1998) 501-511] on the Hopf
invariant one problem and Morin mappings, this implies that a manifold
P^p with odd Euler characteristic does not admit Morin mappings into
R^{2k+1} for p > 2k not equal to 1,3 or 7.}



\maketitle

\section{Introduction}
Let $P$ and $Q$ be two smooth manifolds of dimensions $p$ and $q$
respectively and suppose that $p\ge q$. The singular points of a
smooth mapping $f\co P\to Q$ are the points of the manifold $P$ at
which the rank of the differential $df$ of the mapping $f$ is less
than $q$. There is a natural stratification breaking the singular
set into finitely many strata. We recall that the kernel rank
$kr_x(f)$ of a smooth mapping $f$ at a point $x$ is the rank of
the kernel of $df$ at $x$. At the first stage of the
stratification every stratum is indexed by a non-negative integer
$i_1$ and defined as $$\Sigma^{i_1}(f)=\{\ x\in P\ |\ kr_x(f)=i_1
\}.$$ The further stratification proceeds by induction. Suppose
that the stratum $\Sigma_{n-1}(f)=\Sigma^{i_1,...,i_{n-1}}(f)$ is
defined. Under assumption that $\Sigma_{n-1}(f)$ is a submanifold
of $P$, we consider the restriction $f_{n-1}$ of the mapping $f$
to $\Sigma_{n-1}(f)$ and define
$$\Sigma^{i_1,...,i_n}(f)=\{\ x\in \Sigma_{n-1}(f)\ |\ kr_x(f_{n-1})=i_n \}.$$
Boardman \cite{Bo} proved that every mapping $f$ can be
approximated by a mapping for which every stratum $\Sigma_n(f)$ is
a manifold.

We abbreviate the sequence $(i_1,...,i_n)$ of $n$ non-negative
integers by $I$. We say that a point of the manifold $P$ is an
{\it $I$-singular point} of a mapping $f$ if it belongs to a
singular submanifold $\Sigma^I(f)$. There is a class of in a sense
the simplest singularities, which are called {\it Morin}. Let
$I_1$ denote the sequence $(p-q+1,0)$ and for every integer $k>1$,
the symbol $I_k$ denote the sequence $(p-q+1,1,...,1,0)$ with $k$
non-zero entries. Then Morin singularities are singularities with
symbols $I_k$. A Morin mapping is an {\it $I_k$-mapping} if it has
no singularities of type $I_{k+1}$. For $k=1,2$ and $3$, points
with the symbols $I_k$ are called {\it fold, cusp} and {\it
swallowtail singular points} respectively. In this terminology,
for example, a fold mapping is a mapping which has only fold
singular points.

Given two manifolds $P$ and $Q$, we are interested in finding a
mapping $P\to Q$ that has as simple singularities as possible. Let
$f\co P\to Q$ be an arbitrary general position mapping. For every
symbol $I$, the $\Z_2$-homology class represented by the closure
$\overline{\Sigma^I(f)}$ does not change under general position
homotopy. Therefore the homology class $[\overline{\Sigma^I(f)}]$
gives an obstruction to elimination of $I$-singularities by
homotopy.

In \cite{Ch} Chess showed that if $p-q$ is odd and $k\ge 4$, then the homology
obstruction corresponding to $I_k$-singularities vanishes. Chess conjectured that in
this case every Morin mapping $f$ is homotopic to a mapping without $I_k$-singular
points.

We will show that the statement of the Chess conjecture holds.
Furthermore we will prove a stronger assertion.

\begin{theorem}\label{t1} Let $P$ and $Q$ be two orientable manifolds,
$p-q$ odd. Then the homotopy class of an arbitrary Morin mapping
$f\co P\to Q$ contains a {\it cusp mapping}.
\end{theorem}

\medskip \noindent{\bf Remark}\qua The standard complex projective
plane $\C P^2$ does not admit a fold mapping \cite{Sae} (see also
\cite{AS}, \cite{Sak}). This shows that the homotopy class of $f$
may contain no mappings with only $I_1$-singularities.

\medskip \noindent{\bf Remark}\qua The assumption on the parity of
the number $p-q$ is essential since in the case where $p-q$ is
even homology obstructions may be nontrivial \cite{Ch}.

\medskip \noindent{\bf Remark}\qua We refer to an excellent review
\cite{SS1} for further comments. In particular, see Remark 4.6,
where the authors indicate that Theorem \ref{t1} does not hold for
non-orientable manifolds.

\medskip In \cite{SS} (see also \cite{KS}) Saeki and Sakuma
describe a remarkable relation between the problem of the
existence of certain Morin mappings and the Hopf invariant one
problem. Using this relation the authors show that if the Euler
characteristic of $P$ is odd, $Q$ is almost parallelizable, and
there exists a cusp mapping $f\co P\to Q$, then the dimension of
$Q$ is $1,2,3,4,7$ or $8$.

Note that if the Euler characteristic of $P$ is odd, then the dimension of $P$ is
even. We obtain the following corollary.

\begin{corollary} Suppose the Euler characteristic of $P$ is odd
and the dimension of an almost parallelizable manifold $Q$ is
odd and different from $1,3,7$. Then there exist no Morin mappings from $P$ into $Q$.
\end{corollary}

\section{Jet bundles and suspension bundles}

Let $P$ and $Q$ be two smooth manifolds of dimensions $p$ and $q$ respectively. {\it
A germ at a point $x\in P$} is a mapping from some neighborhood about $x$ in $P$
into $Q$. Two germs are {\it equivalent} if they coincide on some neighborhood of
$x$. The class of equivalence of germs (or simply the germ) at $x$ represented by a
mapping $f$ is denoted by $[f]_x$.

Let $U$ be a neighborhood of $x$ in $P$ and $V$ be a neighborhood
of $y=f(x)$ in $Q$. Let $$\tau_U\co (U,x)\to (\R^p,0) \hspace{5mm}
\mathrm{and} \hspace{5mm} \tau_V\co (V,y)\to (\R^q,0)$$ be
coordinate systems. Two germs $[f]_x$ and $[g]_x$ are {\it
$k$-equivalent} if the mappings $\tau_V\circ f \circ \tau_U^{-1}$
and $\tau_V\circ g\circ\tau_U^{-1}$, which are defined in a
neighborhood of $0\in \R^p$, have the same derivatives at $0\in
\R^p$ of order $\le k$. The notion of $k$-equivalence is
well-defined, i.e. it does not depend on choice of representatives
of germs and on choice of coordinate systems. A class of
$k$-equivalent germs at $x$ is called {\it a $k$-jet}. The set of
all $k$-jets constitute a set $J^k(P,Q)$. The projection
$J^k(P,Q)\to P\times Q$ that takes a germ $[f]_x$ into a point
$x\times f(x)$ turns $J^k(P,Q)$ into a bundle (for details see
\cite{Bo}), which is called {\it the $k$-jet bundle over $P\times
Q$.}

Let $y$ be a point of a manifold and $V$ a neighborhood of $y$. We
say that two functions on $V$ lead to the same local function at
$y$, if at the point $y$ their partial derivatives agree. Thus a
local function is an equivalence class of functions defined on a
neighborhood of $y$. The set of all local functions at the point
$y$ constitutes an algebra of jets $\mathcal F(y)$. Every smooth
mapping $f\co (U,x)\to (V,y)$ defines a homomorphism of algebras
$f^*\co  \mathcal F(y) \to \mathcal F(x).$ The maximal ideal $m_y$
of $\mathcal F(y)$ maps under the homomorphism $f^*$ to the
maximal ideal $m_x\subset\mathcal F(x)$. The restriction of $f^*$
to $m_y$ and the projection of $f^*(m_y)\subset m_x$ onto
$m_x/m_x^{k+1}$ lead to a homomorphism
$$f_{k,x}\co  m_y \to m_x/m_x^{k+1}.$$
It is easy to verify that $k$-jets of mappings $(U,x)\to (V,y)$
are in bijective correspondence with algebra homomorphisms $m_y\to
m_x/m_x^{k+1}$. That is why we will identify a $k$-jet with the
corresponding homomorphism.

The projections of $P\times Q$ onto the factors induce from the
tangent bundles $TP$ and $TQ$ two vector bundles $\xi$ and $\eta$
over $P\times Q$. The latter bundles determine a bundle $\mathcal
{HOM}(\xi,\eta)$ over $P\times Q$. The fiber of
$\mathcal{HOM}(\xi,\eta)$ over a point $x\times y$ is the set of
homomorphisms $Hom(\xi_x,\eta_y)$ between the fibers of the
bundles $\xi$ and $\eta$. The bundle $\xi$ determines the $k$-th
symmetric tensor product bundle $\circ^k\xi$ over $P\times Q$,
which together with $\eta$ leads to a bundle
$\mathcal{HOM}(\circ^k\xi,\eta)$.

\begin{lemma} The $k$-jet bundle contains a vector subbundle $\mathcal C^k$
isomorphic to $\mathcal{HOM}(\circ^k\xi,\eta)$.
\end{lemma}
\begin{proof} Define $\mathcal C^k$ as the union of those $k$-jets
$f_{k,x}$ which take $m_y$ to $m^k_x$. With each $f_{k,x}\in C^k$
we associate a homomorphism  (for details, see \cite[Theorem
4.1]{Bo})
\begin{equation}\label{eq1}
\underbrace{\xi_x\circ...\circ \xi_x}_{\mbox{$k$}}\otimes
m_y/m_y^2\to \R
\end{equation}
which sends $v_1\circ....\circ v_k\otimes\alpha$ into the value of
$v_1\circ...\circ v_k$ at a function representing
$f_{k,x}(\alpha).$ In view of the isomorphism $m_y/m_y^2 \approx
Hom(\eta_y,\R)$, the homomorphism (\ref{eq1}) is an element of
$Hom(\circ^k\xi_x,\eta_y)$. It is easy to verify that the obtained
correspondence $C^k\to\mathcal{HOM}(\circ^k\xi_x,\eta_y)$ is an
isomorphism of vector bundles. \end{proof}

\begin{corollary}\label{c1}
There is an isomorphism $J^{k-1}(P,Q)\oplus \mathcal C^k\approx
J^k(P,Q)$.
\end{corollary}
\begin{proof} Though the sum of two algebra
homomorphisms may not be an algebra homomorphism, the sum of a
homomorphism $f_{k,x}\in J^k(P,Q)$ and a homomorphism $h\in
\mathcal C^k$ is a well defined homomorphism of algebras
$(f_{k,x}+h)\in J^k(P,Q)$. This defines an action of $\mathcal
C^k$ on $J^k(P,Q)$. Two $k$-jets $\alpha$ and $\beta$ map under
the canonical projection
$$J^k(P,Q)\longrightarrow J^{k}(P,Q)/\mathcal C^k$$ onto one point
if and only if $\alpha$ and $\beta$ have the same $(k-1)$-jet.
Therefore $J^k(P,Q)/\mathcal C^k$ is canonically isomorphic to
$J^{k-1}(P,Q)$. \end{proof}

\noindent{\bf Remark}\qua The isomorphism $J^{k-1}(P,Q)\oplus \mathcal C^k\approx
J^{k}(P,Q)$ constructed in Corollary \ref{c1} is not canonical, since there is no
canonical projection of the $k$-jet bundle onto $\mathcal C^k$.

In \cite{Ro} Ronga introduced the bundle
$$S^k(\xi,\eta)=\mathcal{HOM}(\xi,\eta)\oplus \mathcal{HOM}(\xi\circ\xi,\eta)\oplus...\oplus \mathcal{HOM}(\circ^k\xi,\eta),$$
which we will call the $k$-suspension bundle over $P\times Q$.

\begin{corollary} The $k$-jet bundle is isomorphic to the $k$-suspension bundle.
\end{corollary}

\section{Submanifolds of singularities}

There are canonical projections $J^{k+1}(P,Q)\to J^{k}(P,Q)$,
which lead to the infinite dimensional {\it jet bundle} $J(P,Q):
=\underleftarrow{lim}\ J^k(P,Q)$. Let $f\co P\to Q$ be a smooth
mapping. Then at every point $x\times f(x)$ of the manifold
$P\times Q$, the mapping $f$ determines a $k$-jet. The $k$-jets
defined by $f$ lead to a mapping $j^kf$ of $P$ to the $k$-jet
bundle. These mappings agree with projections of
$\underleftarrow{lim}\ J^k(P,Q)$ and therefore define a mapping
$jf\co P\to J(P,Q)$, which is called the jet extension of $f$. We
will call a subset of $J(P,Q)$ {\it a submanifold of the jet
bundle} if it is the inverse image of a submanifold of some
$k$-jet bundle. A function $\Phi$ on the jet bundle is said to be
{\it smooth} if locally $\Phi$ is the composition of the
projection onto some $k$-jet bundle and a smooth function on
$J^k(P,Q)$. In particular, the composition $\Phi\circ jf$ of a
smooth function $\Phi$ on $J(P,Q)$ and a jet extension $jf$ is
smooth. {\it A tangent to the jet bundle vector} is a differential
operator. {\it A tangent to $J(P,Q)$ bundle} is defined as a union
of all vectors tangent to the jet bundle.

Suppose that at a point $x\in P$ the mapping $f$ determines a jet
$z$. Then the differential of $jf$ sends differential operators at
$x$ to differential operators at $z$, that is $d(jf)$ maps $T_xP$
into some space $D_z$ tangent to the jet bundle. In fact, the
space $D_z$ and the isomorphism $T_xP\to D_z$ do not depend on
representative $f$ of the jet $z$. Let $\pi$ denote the
composition of the jet bundle projection and the projection of
$P\times Q$ onto the first factor. Then the tangent bundle of the
jet space contains a subbundle $D$, called {\it the total tangent
bundle}, which can be identified with the induced bundle $\pi^*TP$
by the property: for any vector field $v$ on an open set $U$ of
$P$, any jet extension $jf$ and any smooth function $\Phi$ on
$J(P,Q)$, the section $V$ of $D$ over $\pi^{-1}(U)$ corresponding
to $v$ satisfies the equation $$V\Phi\circ jf=v(\Phi\circ jf).$$
We recall that the projections $P\times Q$ onto the factors induce
two vector bundles $\xi$ and $\eta$ over $P\times Q$ which
determine a bundle $\mathcal {HOM}(\xi,\eta)$. There is a
canonical isomorphism between the $1$-jet bundle and the bundle
$\mathcal{HOM}(\xi,\eta)$. Consequently $1$-jet component of a
$k$-jet $z$ at a point $x\in P$ defines a homomorphism $h\co
T_xP\to T_yQ$, $y=z(x)$. We denote the kernel of the homomorphism
$h$ by $K_{1,z}$. Identifying the space $T_xP$ with the fiber
$D_z$ of $D$, we may assume that $K_{1,z}$ is a subspace of $D_z$.
Hence at every point $z\in J(P,Q)$ we have a space $K_{1,z}$.
Boardman showed that the union $\Sigma^i=\Sigma^i(P,Q)$ of jets
$z$ with $dim\: K_{1,z}=i$ is a submanifold of $J(P,Q)$.

Suppose that we have already defined a submanifold
$\Sigma_{n-1}=\Sigma^{i_1,...,i_{n-1}}$ of the jet space. Suppose
also that at every point $z\in \Sigma_{n-1}$ we have already
defined a space $K_{n-1,z}$. Then the space $K_{n,z}$ is defined
as $K_{n-1,z}\cap T_z\Sigma_{n-1}$ and $\Sigma_n$ is defined as
the set of points $z\in \Sigma_{n-1}$ such that $dim\:
K_{n,z}=i_n$. Boardman proved that the sets $\Sigma_n$ are
submanifolds of $J(P,Q)$. In particular every submanifold
$\Sigma_n$ comes from a submanifold of an appropriate finite
dimensional $k$-jet space. In fact the submanifold with symbol
$I_n$ is the inverse image of the projection of the jet space onto
$n$-jet bundle. To simplify notation, we denote the projections of
$\Sigma_n$ to the $k$-jet bundles with $k\ge n$ by the same symbol
$\Sigma_n$.

Let us now turn to the $k$-suspension bundle. Following the paper
\cite{Bo}, we will define submanifolds $\tilde \Sigma^I$ of the
$k$-suspension bundle.

A point of the $k$-suspension bundle over a point $x\times y\in P
\times Q$ is the set of homomorphisms $h=(h_1,...,h_k)$, where
$h_i\in Hom(\circ^i\xi_x,\eta_y)$. For every $k$-suspension $h$ we
will define a sequence of subspaces $T_xP=K_0\supset
K_1\supset...\supset K_k$. Then we will define the singular set
$\tilde\Sigma^{i_1,...,i_n}$ as
$$\tilde\Sigma^{i_1,...,i_n}=\{\ h\ |\ dim\:K_j=i_j\ \mathrm{for}\ j=1,...,n\ \}.$$
We start with definition of a space $K_1\supset K_0$ and a
projection of $P_0=T_yQ$ onto a factor space $Q_1$. The
$h_1$-component of $h$ is a homomorphism of $K_0$ into $P_0$. We
define $K_1$ and $Q_1$ as the kernel and the cokernel of $h_1$:
$$0\longrightarrow K_1\longrightarrow K_0\stackrel{h_1}\longrightarrow P_0\longrightarrow Q_1\longrightarrow 0.$$
The cokernel homomorphism of this exact sequence gives rise to a homomorphism
$Hom(K_1,P_0)\to Hom(K_1,Q_1)$, coimage of which is denoted by $P_1$. The sequence
of the homomorphisms
$$Hom(K_1\circ K_1,P_0)\to Hom(K_1,Hom(K_1,P_0))\to Hom(K_1,P_1)$$
takes the restriction of $h_2$ on $K_1\circ K_1$ to a homomorphism
$\sigma(h_2)\co K_1\to P_1$. Again the spaces $K_2$ and $Q_2$ are
respectively defined as the kernel and the cokernel of the
homomorphism $\sigma(h_2)$.

The definition continues by induction. In the $n$-th step we are given some spaces
$K_i,Q_i$ for $i\le n$, spaces $P_i$ for $i\le n-1$ and projections
$$Hom(K^{n-1},P_0)\to P_{n-1},$$
$$P_{n-1}\to Q_n,$$
where $K^{n-1}$ abbreviates the product $K_{n-1}\circ...\circ
K_1$.

First we define $P_n$ as the coimage of the composition
$$Hom(K^n,P_0)\to
Hom(K_n,Hom(K^{n-1},P_0))\to Hom(K_n,Q_n),$$ where the latter
homomorphism is determined by the two given projections. Then we
transfer the restriction of the homomorphism $h_{n+1}$ on
$K_n\circ K^n$ to a homomorphism $\sigma(h_{n+1})\co K_n\to P_n$
using the composition
$$Hom(K_n\circ K^n,P_0)\to
Hom(K_n,Hom(K^n,P_0))\to Hom(K_n,P_n).$$ Finally we define $K_{n+1}$ and $Q_{n+1}$
by the exact sequence
$$0\longrightarrow K_{n+1}\longrightarrow K_n\stackrel{\sigma(h_{n+1})}\longrightarrow
P_n\longrightarrow Q_{n+1}\longrightarrow 0.$$ In the previous
section we established a homeomorphism between the fibers of the
$k$-jet bundle and $k$-suspension bundle. Suppose that
neighborhoods of points $x\in P$ and $y\in Q$ are equipped with
coordinate systems. Then every $k$-jet $g$ which takes $x$ to $y$
has the canonical decomposition into the sum of $k$-jets $g_i$,
$i=1,...,k$, such that in the selected coordinates the partial
derivatives of the jet $g_i$ at $x$ of order $\ne i$ and $\le k$
are trivial. In other words the choice of local coordinates
determines a homeomorphism
\begin{equation}\label{e9}
J^k(P,Q)|_{x\times y}\to \mathcal C^1|_{x\times y}\oplus ...
\oplus \mathcal C^k|_{x\times y}.
\end{equation}
Since $\mathcal C^i|_{x\times y}$ is isomorphic to
$Hom(\circ^i\xi_x,\eta_y)$, we obtain a homeomorphism between the
fibers of the $k$-jet bundle and $k$-suspension bundle.

\medskip \noindent{\bf Remark}\qua From \cite{Bo} we deduce that
this homeomorphism takes the singular submanifolds $\Sigma^I$ to
$\tilde\Sigma^I$. Suppose that a $k$-jet $z$ maps onto a
$k$-suspension $h=(h_1,...,h_k)$. The homomorphisms $\{h_i\}$
depends not only on $z$ but also on choice of coordinates in
$U_i$. However Boardman \cite{Bo} showed that the spaces $K_i$,
$Q_i$, $P_i$ and the homomorphisms $\sigma(h_i)$ defined by $h$
are independent from the choice of coordinates.

\begin{lemma}\label{equiv}  For every integer $k\ge 1$, there is a homeomorphism of bundles
$r_k\co  J^k(P,Q)\to S^k(\xi,\eta)$ which takes the singular sets
$\Sigma^I$ to $\tilde\Sigma^I$.
\end{lemma}
\begin{proof} Choose covers of $P$ and $Q$ by closed discs. Let
$U_1,...,U_t$  be the closed discs of the product cover of
$P\times Q$. For each disc $U_i$, choose a coordinate system which
comes from some coordinate systems of the two disc factors of
$U_i$. We will write $J^k$ for the $k$-jet bundle and $J^k|_{U_i}$
for its restriction on $U_i$. We adopt similar notations for the
$k$-suspension bundle. The choice of coordinates in $U_i$ leads to
a homeomorphism
$$\beta_i\co J^k|_{U_i}\to S^k|_{U_i}.$$
Let $\{\varphi_i\}$ be a partition of unity for the cover
$\{U_i\}$ of $P\times Q$. We define $r_k\co J^k\to S^k$ by
$$r_k=\varphi_1\beta_1+\varphi_2\beta_2+...+\varphi_k\beta_k.$$
Suppose that $U_i\cap U_j$ is nonempty and $z$ is a $k$-jet at a
point of $U_i\cap U_j$. Suppose
$$\beta_i(z)=(h^i_1,...,h^i_k)\hspace{3mm}
\mathrm{and}\hspace{3mm} \beta_j(z)=(h^j_1,...,h^j_k).$$ Then by
the remark preceding the lemma, the homomorphisms $\sigma(h^i_s)$
and $\sigma(h^j_s)$ coincide for all $s=1,...,k$. Consequently,
$r_k$ takes $\Sigma^I$ to $\tilde\Sigma^I$.

The mapping $r_k$ is continuous and open. Hence to prove that
$r_k$ is a homeomorphism it suffices to show that $r_k$ is
one-to-one.

For $k=1$, the mapping $r_k$ is the canonical isomorphism. Suppose
that $r_{k-1}$ is one-to-one and for some different $k$-jets $z_1$
and $z_2$, we have $r_k(z_1)=r_k(z_2)$. Since $r_{k-1}$ is
one-to-one, the $k$-jets $z_1$ and $z_2$ have the same $(k-1)$-jet
components. Hence there is $v\in \mathcal C^k$ for which
$z_1=z_2+v$. Here we invoke the fact that $\mathcal C^k$ has a
canonical action on $J^k$.

For every $i$, we have $\beta_i(z_1)=\beta_i(z_2)+\beta_i(v)$.
Therefore
\begin{equation}\label{r}
r_k(z_1)=r_k(z_2)+r_k(v).
\end{equation}
The restriction of the mapping $r_k$ to ${\mathcal C^k}$ is a
canonical identification of $\mathcal C^k$ with
$\mathcal{HOM}(\circ^k\xi_k,\eta)$. Hence $r_k(v)\ne 0$. Then
(\ref{r}) implies that $r_k(z_1)\ne r_k(z_2)$.
\end{proof}

\begin{corollary} There is an isomorphism of bundles $r\co J(P,Q)\to S(\xi,\eta)$ which
takes every set $\Sigma_n$ isomorphically onto $\tilde\Sigma_n$.
\end{corollary}

The space $J^k(P,Q)$ may be also viewed as a bundle over $P$ with projection
$$\pi\co J^k(P,Q)\to P\times Q\to P.$$ Let $f\co P\to Q$ be a smooth mapping. Then at
every point $p\in P$ the mapping $f$ defines a $k$-jet.
Consequently, every mapping $f\co P\to Q$ gives rise to a section
$j^kf\co P\to J^k(P,Q),$ which is called {\it the $k$-extension of
$f$} or {\it the $k$-jet section afforded by $f$}. The sections
$\{j^kf\}_k$ determined by a smooth mapping $f$ commute with the
canonical projections $J^{k+1}(P,Q)\to J^{k}(P,Q)$. Therefore
every smooth mapping $f\co P\to Q$ also defines a section $jf\co
P\to J(P,Q)$, which is called the jet extension of $f$.

A smooth mapping $f$ is {\it in general position} if its jet extension is
transversal to every singular submanifold $\Sigma^I$. By the Thom Theorem every
mapping has a general position approximation.

Let $f$ be a general position mapping. Then the subsets
$(jf)^{-1}(\Sigma^I)$ are submanifolds of $P$. Every condition
$kr_x(f_{n-1})=i_n$ in the definition of $\Sigma^I(f)$ can be
substituted by the equivalent condition $dim\: K_{n,x}(f) = i_n$,
where the space $K_{n,x}(f)$ is the intersection of the kernel of
$df$ at $x$ and the tangent space $T_x\Sigma_{n-1}(f)$. Hence the
sets $(jf)^{-1}(\Sigma^I)$ coincide with the sets $\Sigma^I(f)$.
In particular the jet extension of a mapping $f$ without
$I$-singularities does not intersect the set $\Sigma^I$.

Let $\Omega_r=\Omega_r(P,Q)\subset J(P,Q)$ denote the union of the regular points
and the Morin singular points with indexes of length at most $r$.

\begin{theorem}[Ando-Eliashberg, \cite{An}, \cite{El}]
Let $f\co P^p\to Q^q$,$p\ge q\ge 2$, be a continuous mapping. The
homotopy class of the mapping $f$ contains an $I_r$-mapping, $r\ge
1$, if and only if there is a section of the bundle $\Omega_r$.
\end{theorem}

Note that every general position mapping $f\co P^p\to Q^q$, $q=1$,
is a fold mapping. That is why for $q=1$, Theorem~\ref{t1} holds
and we will assume that $q\ge 2$.

Let $\tilde\Omega_r$ denote the subset of the suspension bundle
corresponding to the set $\Omega_r(P,Q)\subset J(P,Q)$. Every
mapping $f\co P\to Q$ defines a section $jf$ of $J(P,Q)$. The
composition $r\circ (jf)$ is a section of $S(P,Q)$. In view of
Lemma \ref{equiv} the Ando-Eliashberg Theorem implies that to
prove that the homotopy class of a mapping $f$ contains a cusp
mapping, it suffices to show that the section of the suspension
bundle defined by $f$ is homotopic to a section of the bundle
$\tilde\Omega_2\subset S(\xi,\eta)$.

\section{Proof of Theorem \ref{t1}}

We recall that in a neighborhood of a fold singular point $x$, the mapping $f$ has
the form
\begin{eqnarray}\label{e1}
T_i&=&t_i,\ \ i=1, 2, ... , q-1, \\
Z  &=&Q(x),\ \ Q(x)=\pm k_1^2\pm...\pm k_{p-q+1}^2. \nonumber 
\nonumber
\end{eqnarray}
If $x$ is an $I_r$-singular point of $f$ and $r>1$, then in some
neighborhood about $x$ the mapping $f$ has the form
\begin{eqnarray}\label{e2}
T_i&=&t_i,\ \ i=1, 2, ... , q-r, \nonumber \\
L_i&=&l_i,\ \ i=2, 3, ... , r,    \\
Z  &=&Q(x) + \sum_{t=2}^{r} l_t k^{t-1} + k^{r+1},\ \ Q(x)=\pm k_1^2\pm...\pm k_{p-q}^2. \nonumber 
\nonumber
\end{eqnarray}
Let $f\co P\to Q$ be a Morin mapping, for which the set
$\Sigma_2(f)$ is nonempty. We define the section $f_i\co P\to
Hom(\circ^i\xi,\eta)$ as the $i$-th component of the section
$r\circ (jf)$ of the suspension bundle $S(\xi,\eta)\to P$. Over
$\overline{\Sigma_2(f)}$ the components $f_1$ and $f_2$ defined by
the mapping $f$ determine the bundles $K_i, Q_i$, $i=1,2$ and the
exact sequences
$$0\longrightarrow K_1\longrightarrow TP\longrightarrow TQ\longrightarrow Q_1\longrightarrow 0,$$
$$0\longrightarrow K_2\longrightarrow K_1\longrightarrow \mathcal{HOM}(K_1,Q_1)
\longrightarrow Q_2\longrightarrow 0.$$ From the latter sequence
one can deduce that the bundle $Q_2$ is canonically isomorphic to
$\mathcal {HOM}(K_2,Q_1)$ and that the homomorphism
\begin{equation}\label{factor}
K_1/K_2\otimes K_1/K_2\longrightarrow Q_1,
\end{equation}
which is defined by the middle homomorphism of the second exact sequence, is a
non-degenerate quadratic form (see Chess, \cite{Ch}). Since the dimension of
$K_1/K_2$ is odd, the quadratic form (\ref{factor}) determines a canonical
orientation of the bundle $Q_1$. In particular the $1$-dimensional bundle $Q_1$ is
trivial. This observation also belongs to Chess \cite{Ch}.

Assume that the bundle $K_2$ is trivial. Then the bundle $Q_2$ being isomorphic to
$\mathcal {HOM}(K_2,Q_1)$ is trivial as well. Let
$$\tilde h\co K_2\to
\mathcal{HOM}(K_2,Q_2)\approx \mathcal{HOM}(K_2\otimes K_2, Q_1)$$
be an isomorphism over $\overline{\Sigma_2(f)}$ and $h\co P\to
\mathcal {HOM}(\circ^3\xi,\eta)$ an arbitrary section, the
restriction of which on $\circ^3 K_2$ over
$\overline{\Sigma_2(f)}$ followed by the projection given by
$\eta\to Q_1$, induces the homomorphism $\tilde h$. Then the
section of a suspension bundle whose first three components are
$f_1,f_2$ and $h$ is a section of the bundle $\tilde\Omega_2$.
Since for $i>0$ the bundle $\mathcal {HOM}(\circ^i\xi,\eta)$ is a
vector bundle, we have that the composition $r\circ(jf)$ is
homotopic to the section $s$ and therefore the original mapping
$f$ is homotopic to a cusp mapping.

Now let us prove the assumption that $K_2$ is trivial over $\overline{\Sigma_2(f)}$.

\begin{lemma}\label{s2} The submanifold
$\overline{\Sigma_2(f)}$ is canonically cooriented in the submanifold
$\overline{\Sigma_1(f)}$.
\end{lemma}
\begin{proof} For non-degenerate quadratic forms of order $n$, we
adopt the convention to identify the index $\lambda$ with the
index $n-\lambda$. Then the index $ind\:Q(x)$ of the quadratic
form $Q(x)$ in (\ref{e1}) and (\ref{e2}) does not depend on choice
of coordinates.

With every $I_k$-singular point $x$ by (\ref{e1}) and (\ref{e2})
we associate a quadratic mapping of the form $Q(x)$. It is easily
verified that for every cusp singular point $y$ and a fold
singular point $x$ of a small neighborhood of $y$, we have
$Q(x)=Q(y)\pm k_{p-q+1}^2$. Moreover, if $x_1$ and $x_2$ are two
fold singular points and there is a path joining $x_1$ with $x_2$
which intersects $\overline{\Sigma_2(f)}$ transversally and at
exactly one point, then $ind\:Q(x_1)-ind\:Q(x_2)=\pm 1$. In
particular, the normal bundle of $\overline{\Sigma_2(f)}$ in
$\overline{\Sigma_1(f)}$ has a canonical orientation.
\end{proof}

\begin{lemma}\label{k2} Over every connected component of $\Sigma_2(f)$ the bundle $K_2$ has a
canonical orientation.
\end{lemma}
\begin{proof} At every point $x\in \overline{\Sigma_2(f)}$ there is
an exact sequence
$$0\longrightarrow K_{3,x}\longrightarrow K_{2,x}\longrightarrow \mathcal{HOM}(K_{2,x},Q_{2,x})
\longrightarrow Q_{3,x}\longrightarrow 0.$$ If the point $x$ is in fact a cusp
singular point, then the space $K_{3,x}$ is trivial and therefore the sequence
reduces to
$$0\longrightarrow K_{2,x}\longrightarrow \mathcal{HOM}(K_{2,x},Q_{2,x})\longrightarrow 0$$
and gives rise to a quadratic form
$$K_{2,x}\otimes K_{2,x}\longrightarrow Q_{2,x}\approx \mathcal{HOM}(K_{2,x},Q_{1,x}).$$
This form being non-degenerate orients the space
$\mathcal{HOM}(K_{2,x},Q_{1,x})$. Since $Q_{1,x}$ has a canonical
orientation, we obtain a canonical orientation of $K_{2,x}$.
\end{proof}

Let $\gamma\co [-1,1]\to \overline{\Sigma_2(f)}$ be a path which
intersects the submanifold of non-cusp singular points
transversally and at exactly one point.
\begin{lemma}\label{change} The canonical orientations of $K_2$ at $\gamma(-1)$ and $\gamma(1)$
lead to different orientations of the trivial bundle $\gamma^*K_2$.
\end{lemma}
\begin{proof} If necessary we slightly modify the path $\gamma$ so
that the unique intersection point of $\gamma$ and the set
$\overline{\Sigma_3(f)}$ is a swallowtail singular point. Then the
statement of the lemma is easily verified using the formulas
(\ref{e2}). \end{proof}

Now we are in position to prove the assumption.

\begin{lemma}\label{l4} The bundle $K_2$ is trivial over $\overline{\Sigma_2(f)}$.
\end{lemma}
\begin{proof} Assume that the statement of the lemma is wrong. Then there
is a closed path $\gamma\co S^1\to \overline{\Sigma_2(f)}$ which
induces a non-orientable bundle $\gamma^*K_2$ over the circle
$S^1$.

We may assume that the path $\gamma$ intersects the submanifold
$\overline{\Sigma_3(f)}$ transversally. Let $t_1,...,t_k,
t_{k+1}=t_1$ be the points of the intersection $\gamma\cap
\overline{\Sigma_3(f)}$. Over every interval $(t_i,t_{i+1})$ the
normal bundle of $\overline{\Sigma_2(f)}$ in
$\overline{\Sigma_1(f)}$ has two orientations. One orientation is
given by Lemma \ref{s2} and another is given by the canonical
orientation of the bundle $K_2$. By Lemma \ref{change} if these
orientations coincide over $(t_{i-1},t_i)$, then they differ over
$(t_i,t_{i+1})$. Therefore the number of the intersection points
is even and the bundle $\gamma^*K_2$ is trivial. Contradiction.
\end{proof}

\noindent{\bf Remark}\qua The statement similar to the assertion of
Lemma~\ref{l4} for the jet bundle $J(P,Q)$ is not correct. The
vector bundle $K_2$ over $\overline{\Sigma^{I_2}}\subset J(P,Q)$
is non-orientable. This follows for example from the study of
topological properties of $\Sigma^{I_r}$ in \cite[\S 4]{An}.

\Addresses\recd


\begin{thebibliography}{99}
\bibitem{AS} {\bf P.~Akhmetev, R.~Sadykov}, {\it A remark on
elimination of singularities for mappings of $4$-manifold into
$3$-manifold}, Top. Appl., 131 (2003), 51-55.
\bibitem{An} {\bf Y.~Ando}, {\it On the elimination of Morin singularities},
J. Math. Soc. Japan, 37 (1985), 471-487; Erratum 39 (1987), 537.
\bibitem{It} {\bf V.~I.~Arnol'd, V.~A.~Vasil'ev, V.~V.~Goryunov, O.~V.~Lyashenko},
{\it Dynamical systems VI. Singularities, local and global
theory}, Encyclopedia of Mathematical Sciences - Vol. 6 (Springer,
Berlin, 1993).
\bibitem{Bo}  {\bf J.~M.~Boardman}, {\it Singularities of differentiable
maps}, Publ. Math., 33 (1967), 21-57.
\bibitem{Ch} {\bf D.~S.~Chess}, {\it A note on the classes $[S_1^k(f)]$}, Proc. Symp. Pure Math., 40 (1983), 221-224.
\bibitem{El} {\bf J.~M.~Eliashberg}, {\it Surgery of singularities of smooth
mappings}, Math. USSR Izv., 6 (1972), 1302-1326.
\bibitem{KS} {\bf S.~Kikuchi, O.~Saeki},
{\it Remarks on the topology of folds}, Proc. Amer. Math. Soc.,
123 (1995), 905-908.
\bibitem{Ro} {\bf F.~Ronga},
{\it Le calculus des classes duales singularit\'{e}s de Boardman
d'ordre deux}, Comment. Math. Helv., 47 (1972), 15-35.
\bibitem{Sae} {\bf O.~Saeki}, {\it Notes on the topology of folds},
J. Math. Soc. Japan, v.44, 3 (1992), 551-566.
\bibitem{SS} {\bf O.~Saeki, K.~Sakuma},
{\it Maps with only Morin singularities and the Hopf invariant one
problem}, Math. Proc. Camb. Phil. Soc., 124 (1998), 501-511.
\bibitem{SS1} {\bf O.~Saeki, K.~Sakuma},
{\it Elimination of Singularities: Thom Polynomials and Beyond},
London Math. Soc., Lecture Notes Ser. 263.
\bibitem{Sak} {\bf K.~Sakuma}, {\it A note on nonremovable cusp
singularities}, Hiroshima Math. J., 31 (2001), 461-465.
\end{thebibliography}
\end{document}